\documentclass{amsart}
\usepackage{graphicx}
\usepackage{color}
\usepackage{amsfonts}
\newtheorem{theorem}{Theorem}
 \newtheorem{definition}{Definition}[section]

\UseRawInputEncoding
\usepackage[russian, english]{babel}

\def\dfrac#1#2{\displaystyle{#1\over #2}}

\def\bV{{\bf V}}

\def\Div{\mbox{div}\,}
\def\Rot{\mbox{rot}\,}

\def\bB{{\bf B}}

\def\bE{{\bf E}}

\begin{document}

\markboth{Rozanova}{The Riemann problem for equations of a cold plasma}

%
%

\title{
The Riemann problem for equations of a cold plasma}

\author{Olga S. Rozanova*}

\address{ Mathematics and Mechanics Department, Lomonosov Moscow State University, Leninskie Gory,
Moscow, 119991,
Russian Federation,
rozanova@mech.math.msu.su}


\subjclass{Primary 35Q60; Secondary 35L60, 35L67, 34M10}
\keywords{Quasilinear hyperbolic system, Riemann problem, non-uniqueness, singular shock,
plasma oscillations}

\maketitle


\begin{abstract} A solution of the Riemann problem is constructed for a nonstrictly hyperbolic inhomogeneous system of equations describing one-dimensional cold plasma oscillations. Each oscillation period includes one rarefaction wave and one shock wave containing a delta singularity. The rarefaction wave can be constructed in a non-unique way, the admissibility principle is proposed.
\end{abstract}
\maketitle


\section{Introduction}

In vector form, the system of hydrodynamic of electron liquid, together with Maxwell's equations, has the form:
\begin{equation}
\label{base1}
\begin{array}{l}
 n_t + \Div(n \bV)=0\,,\quad
\bV_ t + \left( \bV \cdot \nabla \right) \bV
=\dfrac {e}{m} \, \left( \bE + \dfrac{1}{c} \left[\bV \times  \bB\right]\right),\vspace{0.5em}\\
\dfrac1{c}  \bE_t = - \dfrac{4 \pi}{c} e n \bV
 + {\rm rot}\, \bB\,,\quad
\dfrac1{c}  \bB_ t  =
 - {\rm rot}\, \bE\,, \quad \Div \bB=0\,,
\end{array}
\end{equation}
where $e, m$ are the charge and mass of the electron (here the electron charge has a negative sign: $ e <0 $),
$ c $ is the speed of light;
$ n,  \bV $ are the density and velocity of
electrons;
$ \bE, \bB $ are the vectors of electric and magnetic fields, $x\in{\mathbb R}^3,$ $t\ge 0$, $\nabla$, $\rm div$, $\Rot$ are the gradient, divergence and vorticity with respect to the spatial variables.
The system of equations  \eqref {base1} is one of the simplest models of plasma, which is often called the equations of hydrodynamics of "cold" plasma, it is well known and described in sufficient detail in textbooks and monographs (see, for example,  \cite {ABR}, \cite {david72}).

 This system has an important subclass of solutions, dependent only on one space variable $x$,
  for which $\bV=(V,0,0)$, $\bE=(E,0,0)$, $\bB\equiv 0$, e.g. \cite{Ch_book}. In dimensionless form it can be written as
 \begin{equation}
\begin{array}{c}
n_t +
\left(n\, V \right)_x
=0,\quad
V_t +
V  V_x =  - E, \quad
E_t = n\, V.
\end{array}
\label{3gl3}
\end{equation}
Assume that the solution is smooth. Then the first and last equations (\ref{3gl3}) imply
$
\left(n +
 E_x \right)_t = 0.
$
For the background density $ n \equiv 1 $ we get
\begin{equation}
 n = 1 -E_x.
\label{Kn}
\end{equation}
This allows us to obtain a hyperbolic system for the two components of the velocity $V$ and the electric field $E$ in the form
\begin{equation}\label{K1}
V_t+VV_x=-E, \quad E_t+VE_x=V,
\end{equation}
where $(V,E)=(V(t,x), E(t,x))$, $t\in {\mathbb R}_+ $,  $x\in {\mathbb R} $. The density $n(t,x)>0$ can be found from \eqref{Kn}.

System \eqref{K1}, \eqref{Kn} can be also rewritten as a pressureless repulsive Euler-Poisson system \cite{ELT}
\begin{eqnarray}\label{EP}
n_t  + (n V)_x=0,\quad
V_t+VV_x =\,  \nabla \Phi, \quad \Delta \Phi =n-n_0,\quad n_0=1,
\end{eqnarray}
where $\Phi$ is a repulsive  force potential, $\nabla \Phi = -E$.

For \eqref{K1} we consider the Cauchy problem
\begin{equation}\label{K2}
(V,E)|_{t=0}=(V_0(x), E_0(x)).
\end{equation}
If the initial data are  $ C ^ 1 $ - smooth functions, then locally in $t$ there exists a smooth solution of \eqref{K1}, \eqref{K2}. Nevertheless, it is known that the derivatives of the solution of such a Cauchy problem can go to infinity for a finite time, which corresponds to the formation of a shock wave, the criterion for the formation of a singularity is known \cite{RCh2021}.
Thus, it makes sense to consider piecewise-smooth functions as the initial data \eqref {K2}, the simplest example of which is the Riemann initial data
\begin{equation}\label{K3}
(V,E)|_{t=0}=(V^0_-+[V]^0 \Theta (x),\,\, E^0_-  +[E]^0 \Theta(x)),
\end{equation}
 where $\Theta (x)$ is the Heaviside function, constants $(V_-, E_-)$ are the values to the left of the jump,  $([V], [E])$  the values to the  jumps, $(V_+=V_-+[V],\,\, E_+=E_-+[E])$ are the values to the right of the jump,  $(V^0_\pm, E^0_\pm)$, $([V]^0, [E]^0)$ are the corresponding values at time zero.
In this case, the density at the initial moment of time is
\begin{equation}\label{n0}
n|_{t=0}=1 -[E]^0 \delta(x).
\end{equation}
Since the initial data contain a delta function, the Riemann problem for the components of the solution $ (V, E, n) $ is singular and the Rankine-Hugoniot conditions cannot be written in the traditional form \cite{Shelkv}.
In order to ensure that the density is positive initially, it is necessary to impose the condition $[E]^0\le 0$.

To construct the shock, we write  system \eqref {K1} in the divergent form
\begin{equation}\label{K4}
n_t + (V n)_x=0, \qquad
\left(\frac{nV^2}{2}+\frac{E^2}{2}\right)_t +\left(\frac{nV^3}{2}\right)_x =0,
\end{equation}
corresponding to the laws of conservation of mass and total energy (for example, \cite {FrCh}).
System \eqref {K4} (together with \eqref{Kn}) is equivalent to \eqref {K1}, \eqref {Kn} for smooth solutions.

The Riemann problem \eqref {K4}, \eqref{Kn}, \eqref {K3}, \eqref {n0} is completely non-standard and demonstrates new phenomena in the construction of both a rarefaction wave and a shock wave.

The difficulty in constructing a solution is associated, in particular, with the fact that  system \eqref {K1} is inhomogeneous and does not have a constant stationary state. To the left and right side of the discontinuity, the solution is a $ 2 \pi $ - periodic function of time. This leads to the fact that the rarefaction wave and the shock wave periodically replace each other. Further,  system \eqref {K1} is hyperbolic, but not strictly hyperbolic: it has the form
\begin{equation*}
u_t + A (u) u_x = f (u),\quad   u = (V, E), \quad  f = (- E , V),
\end{equation*}
the matrix $ A $ has a complete set of eigenvectors with coinciding eigenvalues $\lambda_1=\lambda_2=V$. Because of this, it has a subclass of solutions in the form of simple waves, distinguished by the condition
 \begin{equation*}\label{KE}
   V^2+E^2= C^2
 \end{equation*}
with a given constant $C$.
We show that this leads to the non-uniqueness of the rarefaction wave for the Riemann problem. Therefore, the question arises about the principles by which one can single out the "correct" solution. In our work, the УcorrectФ one is chosen for which the total energy density is minimal. 

 When constructing a singular shock wave, we use  homogeneous conservative system of two equations \eqref {K4},  which are linked by the differential relation \eqref {Kn}. This formulation has not been encountered before, although a modification of the method previously used for the case of equations of the pressureless gas dynamics with energy \cite{NRS} can be used to construct a solution to the Riemann problem. The shock wave satisfies the so-called ``supercompression" conditions, which are traditionally used to distinguish admissible singular shock waves \cite{Shelkv}.

The paper is organized as follows. In Sec.\ref{Char} we discuss the structure of characteristics which is crucial for construction of rarefaction and shock waves. In Sec.\ref{Rar} we construct the rarefaction wave for Riemann data \eqref{K3} of a general form and then show that for the data corresponding to a simple wave the rarefaction can be constructed non-uniquely. We also propose two variational conditions of admissibility of the rarefaction waves for this case.  In Sec.\ref{Shock} we give a definition of the strongly singular solution for an arbitrary piecewise smooth initial data, prove an analog of the Rankine-Hugoniot conditions (Theorem 1), study the mass and energy transfer for a singular shock wave (Theorem 2). Then we construct the singular shock for  piecewise smooth initial data \eqref {K3} and give two examples. The first example corresponds to the case of simple wave, here we compare the result obtained starting from conservative form
\eqref {K4} and the result obtained from a divergence form, natural for the Hopf equation. The second example show how it is possible to
construct the shock in the case where the shock position has an extremum on the characteristic plane.  Sec.\ref{Dis} contains a discussion about a physical and mathematical sense of the results obtained and mention works concerning shock waves in plasma for other models.

 \section{Characteristics}\label{Char}

\begin{center}
\begin{figure}[htb]
\includegraphics[scale=0.6]{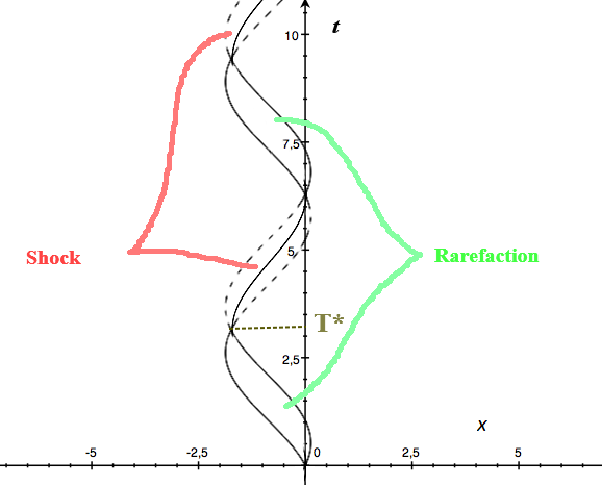}
\caption{ Characteristics and their intersections: rarefaction waves and shock waves. }\label{Pic1}
\end{figure}
\end{center}

 The equations for the characteristics corresponding to system \eqref {K1}  have the form
 \begin{equation}\label{Kchar}
 \dfrac {d V}{dt}= -E,\quad \dfrac {d E}{dt}= V,\quad \dfrac {d x}{dt}= V,
 \end{equation}
whence, first, it follows that along the characteristics
 \begin{equation}\label{KEloc}
   \dfrac {d (V^2+E^2)}{dt}= 0,
 \end{equation}
 and also, according to \eqref{K3},
  \begin{eqnarray*}\label{KvE}
   V_\pm(t)&=-E^0_\pm \sin t + V^0_\pm \cos t, \quad & E_\pm(t)=V^0_\pm \sin t + E^0_\pm \cos t,\\
    x_\pm(t)&=V^0_\pm \sin t + E^0_\pm ( \cos t -1)+ x_0,&\quad x_0=0.\label{Kx}
 \end{eqnarray*}
 It is easy to see that for $ [E]^0 \ne 0 $ the characteristics $ x_- (t) $ and $ x_+ (t) $, corresponding to the states to the left and to the right of the discontinuity, intersect once inside each period $ 2 \pi $.  Therefore, on that part of the period where $ x_- (t) <x_+ (t) $, it is necessary to construct a continuous solution, and on the part where $ x_- (t)> x_+ (t) $, that is, there is an intersection of the characteristics, we construct a shock wave. The moment of time at which $ x_- (t) = x_ + (t) $, we denote by $ T_* $, $ T_* \in (0, 2 \pi) $. Fig.\ref{Pic1}~gives a schematic representation of the behavior of the characteristics, where the rarefaction wave comes first.

Note that \eqref {KEloc} implies that the value $ C^2 = V^2 + E^2 $ is constant for each specific characteristic, but in general it is a function of $ t $ and $ x$.

\section{Construction of a rarefaction wave}\label{Rar}

Suppose the initial data is such that $ V^0_- <V^0_+ $, that is, $ x_- (t) <x_+ (t) $,  and first the initial data \eqref {K3} generate a rarefaction wave.
Between the characteristics $ x_-(t) $ and $ x_+(t) $, it is necessary to construct a continuous solution $ (V, E) $ connecting the states $ (V_- (t), E_- (t)) $ and $ (V_+ (t), E_+ (t)) $. Recall that the moment of time at which $ x_- (t) = x_ + (t) $, we denote by $ T_* $, $ T_* \in (0, 2 \pi) $.

The rarefaction wave, of course, is not a smooth solution, it satisfies the conservative system \eqref {K4} with the additional condition \eqref{Kn} in the usual sense of the integral identity.

\subsection{The linear profile solution}
It is easy to check that a continuous solution $(V,E)$ can be constructed by joining the states $(V_-(t), E_-(t))$ and $(V_+(t), E_+( t))$, between characteristics with the help of functions linear in $x$, i.e.
\begin{equation}\label{K7}
(V,E)=\left\{\begin{array}{ll} (V_-(t),E_-(t)), & x<x_-(t);\\
(V_{r_1},E_{r_1})= (a(t) x + b(t), c(t) x + d(t)),& x\in [x_-(t),x_+(t)];\\
(V_+(t),E_+(t)), &  x>x_+(t),
\end{array}\right.
\end{equation}
 with
 \begin{equation}\label{K8}
 a(t)=\frac{[E]^0\sin t -[V]^0\cos t}{-[V]^0\sin t +[E]^0(1-\cos t)}, \quad c(t)=\frac{-[V]^0\sin t -[E]^0\cos t}{-[V]^0\sin t +[E]^0(1-\cos t)},
 \end{equation}
\begin{equation}\label{K9}
 b(t)=\frac{(V_+^0 E_-^0 -E_+^0 V_-^0)(1-\cos t)}{-[V]^0\sin t +[E]^0(1-\cos t)}, \quad d(t)=\frac{(E_+^0 V_-^0-V_+^0 E_-^0 )\sin t}{-[V]^0\sin t +[E]^0(1-\cos t)}.
 \end{equation}
 Then
\begin{equation*}
n = 1 - c (t) \chi_ {(x_- (t), x_+ (t))},
  \end{equation*}
 where $ \chi_{(x_- (t), x_+ (t))} $ is the characteristic function of the interval $ (x_- (t), x_+ (t)) $, for $ t \in (0, T_*) $ the density does not contain a delta function, but the singular component that was present in the initial data is again formed at $ t = T_* $.

\subsection{Simple waves}


  The system \eqref {K1} has a subclass of solutions distinguished by the condition
 \begin{equation}\label{KE}
V^2+E^2= C^2 (\equiv \rm const)
 \end{equation}
with a given constant $ C $, the so called simple waves. In this case, \eqref {K1} reduces on smooth solutions to one equation
\begin{equation}\label{K5}
V_t+V V_x=-\sigma\sqrt{C^2-V^2},\quad \sigma={\rm sign} (-V_x)=\pm 1,\quad E=\sigma\sqrt{C^2-V^2},
\end{equation}
moreover, $ V_{xx} \ne 0 $ on no set of positive measure. The last requirement means that the solution cannot become constant on any interval, but at the points at which $ C^2 = V^2 $ the value of $ \sigma $ changes its sign to the opposite. The second conservation law \eqref {K4} in this situation turns out to be a consequence of the first.

In the initial conditions~\eqref{K3} the values $ E^0_- $ and $ E^0_+ $ are expressed as $ E^0_- = \pm \sqrt {C^2- (V^0_-)^2} $, $ E_+^0 = \pm \sqrt {C^2- (V_+^0)^2} $,
so as to ensure the condition $ [E]^0 \le 0 $.

 It is easy to see that a function of the form \eqref {K7} with an intermediate state $ (V_ {r_1}, E_ {r_1}) $ is not a solution to the equation \eqref {K5}. Let us show that in this case another continuous solution can be constructed, with another function $ (V_ {r_2}, E_ {r_2}) $ as an intermediate state.

  Indeed, the general solution of \eqref {K5}, written implicitly, looks like
 \begin{equation*}\label{K10}
   x-\sigma\sqrt{C^2-V^2}=F\left(t+\arctan \frac{V}{\sigma\sqrt{C^2-V^2}}\right),
 \end{equation*}
 with an arbitrary smooth function $ F $. In order to find the function $ F $ corresponding to the initial data \eqref {K3}, \eqref {KE}, we will construct the function $ X (t, V) $ inverse to $ V_{r_2} (t, x) $ for every fixed $ t \in (0, T _ *) $. For $ t = 0 $ such a function is multivalued.

We require that for $ t = 0 $ the condition $ X (0, V) = 0 $ holds for $ V \in (V_-^0, V_+^0) $. Then $ F = - \tan \frac {\sqrt {C ^ 2- \xi ^2}} {\xi} $, $ \xi = \sigma \sqrt {C^2-V^2} $. After transformations, we get
 \begin{eqnarray}\label{K11_1}
 &&X_1(t,V)=C(\cos q -\cos (q+t)), \, \mbox{если} \, (V_\pm)_t<0\quad (\sigma=1),\\
 && X_2(t,V)=C(-\cos q +\cos (q-t)), \,\mbox{если} \,(V_\pm)_t>0\quad (\sigma=-1), \label{K11_2}\\
  && q=\arcsin\frac{V}{C}.\nonumber
 \end{eqnarray}
Note that in each case the monotonicity of $ V $ in $ x $ ensures the existence of an inverse function.

The situation is considered separately when the behavior of the solution between the right and left characteristics is given by different formulas. Namely, consider the time $ T_1 $ at which $ V_-'(t) = 0 $ and the time $ T_2 $ at which $ V_+'(t) = 0 $. Between $ T_1 $ and $ T_2 $ there is a moment $ T_0 $, at which $ V_+ (t) = V_- (t) $, and therefore, the jump disappears. However, at such a point the characteristics do not intersect, that is, $ x_+ (T_0) \ne x_- (T_0) $. To construct a continuous solution in such a situation, we need auxiliary curves $ X_1 (t, q_-) $ and $ X_2 (t, q_+) $, where $ q_{\pm} = \arcsin \frac {V^0_{\pm}}{C} $.

Then for $ t \in (0, T_*) $ the continuous solution of  problem  \eqref {K3}, \eqref {K5}, \eqref {KE} can be written as
\begin{equation}\label{K13}
  V_{r_2}(t,x)=\left\{\begin{array}{cc}V_-(t), &x<X_-(t),\\
  V_i(t,x), &X_-(t)<x<X_+(t),\\V_+(t), &x>X_+(t), \end{array}\right.
\end{equation}
where
$$X_-(t)=\left\{\begin{array}{cc}x_-(t), &t<T_1, \, t>T_0\\
 X_2(t,q_+), & t\in [T_1,T_0] \end{array}\right., $$
 $$X_+(t)=\left\{\begin{array}{cc}x_+(t), &t<T_0, \, t>T_2\\
 X_1(t,q_-), & t\in [T_0,T_2] \end{array}\right., $$ and $  V_i(t,x)$ is the function inverse to $ X_i (t, V), $ $ i = 1,2 $, given by formulas \eqref{K11_1}, \eqref{K11_2}.

\begin{center}

\begin{figure}[htb]

\hspace{-1cm}
\hspace{1cm}
\begin{minipage}{0.4\columnwidth}
\includegraphics[scale=0.7]{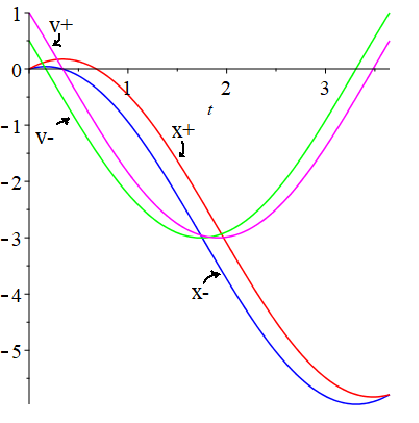}
\end{minipage}
\hspace{1.5cm}
\begin{minipage}{0.4\columnwidth}
\includegraphics[scale=0.35]{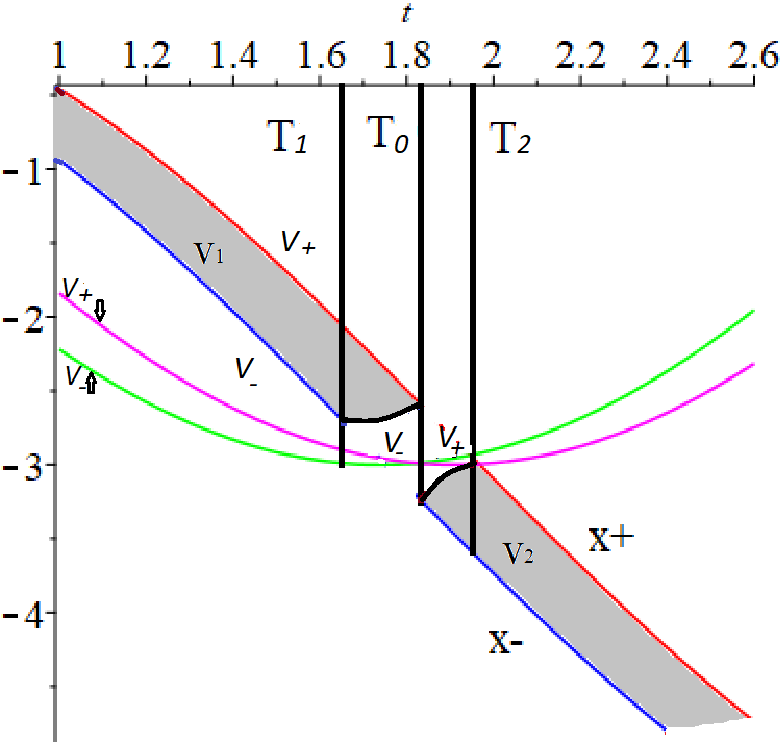}
\end{minipage}
\caption{ Characteristics and and values $V_-$ and $V_+$ from the left and right side of the rarefaction wave. }\label{Pic2}
\end{figure}
\end{center}

Thus, a continuous solution to the problem \eqref {K1}, \eqref {K3}, \eqref {KE} can be constructed as
\begin{equation}\label{K15a}
 (V, E)=\left\{\begin{array}{cc}(V_-(t), E_-(t)), &x<x_-(t),\\
  (V_{r_2}, E_{r_2}), &x_-(t)<x<x_+(t),\\(V_+(t), E_-(t)), &x>x_+(t), \end{array}\right.
\end{equation}
 where $ (V_{r_2}, E_{r_2}) $, where $ V_{r_2} $ is given by \eqref {K13}, and $ E_{r_2} = \pm \sqrt {C^2-V^2_{r_2 }} $, the sign matches the one that was selected in the initial data \eqref {K3}. Fig.\ref{Pic2} presents the construction of the rarefaction wave on the characteristic plane.

\subsection{Nonuniqueness of rarefaction wave}

 Obviously, \eqref {K7} and \eqref{K15a} are different continuous solutions. Moreover, on their basis it is possible to construct an infinite number of other rarefaction waves. Indeed, one can check that $ V_{r_2} $ is an upward convex function and, for $ t = t_1 \in (0, T_*) $, we can choose any point $ x_1 \in (x_- (t_1), x_+ (t_1 )) $ and replace on the segment $ (x_- (t_1), x_1) $ by a linear function. Next, we find the position of the right point of the linear segment as a solution to the problem for $ t \in (0, T_*) $ as $ \dot x = V_ {r_2} (t, x) $, $ x (t_1) = x_1 $. Such linear sections can be built in any number.

\subsection{Admissibility of the rarefaction wave}

The question of choosing the ``correct'' continuous solution can be solved proceeding from the minimality of the total energy of the rarefaction wave
\begin{equation*}\label{local_E}
 {\mathbb E}(t) =\frac12 \,\int\limits_{x_-(t)}^{x_+(t)} (n V^2 + E^2) \, dx,
\end{equation*}
 see \eqref{K4}.

  For the solution $ (V_{r_2}, E_{r_2}) $
   \begin{eqnarray*}\label{local_E}
 {\mathbb E}(t) =\frac12 \,\int\limits_{x_-(t)}^{x_+(t)} ((1-E_x) (C^2- E^2) + E^2) \, dx=
\frac12(C^2 \Delta x - C^2 [E] + \frac13 [E^3]),
\end{eqnarray*}
   where $\Delta x={x_+(t)}-{x_-(t)}\ge 0$, $[E]=E_+-E_-=\Delta x + [E]^0$, $[E^3]=(E_+)^3-(E_-)^3$.

   For the solution $ (V_{r_1}, E_{r_1}) $
   \begin{eqnarray*}\label{local_E}
&& {\mathbb E}(t) =\frac12 \,\int\limits_{x_-(t)}^{x_+(t)} ((1-c) (a x +b)^2 + (c x +d)^2) \, dx, 
\end{eqnarray*}
 where $a, b, c, d$ given as \eqref{K8}, \eqref{K9}.

 It can be readily computed that
  \begin{eqnarray*}\label{local_E}
&& {\mathbb E}(V_{r_2}, E_{r_2})-{\mathbb E}(V_{r_1}, E_{r_1}) =-\frac16 [E]^0 (C^2-(E^0_+E^0_-+V^0_+V^0_-) =\\
&&-\frac{1}{12} [E]^0 (([E]^0)^2+([V]^0)^2)  \ge 0,   \quad t\in (0, T_*).
\end{eqnarray*}
Here we take into account $[E]^0\le 0$ and $(E_+)^2+(V_+)^2=(E_-)^2+(V_-)^2=C^2$.

Thus, if $[E]^0< 0$,  for reasons of less energy ${\mathbb E}$ we have to choose $(V_{r_1}, E_{r_1})$.

\bigskip

 2. Another way to distinguish an acceptable rarefaction wave is the pointwise minimality of the local energy
  \begin{equation*}\label{local_E}
 {\mathcal E}(t,x) = V^2 + E^2.
\end{equation*}
Indeed, $${\mathcal E}(V_{r_2}, E_{r_2})= V^2_{r_2}+ E^2_{r_2}= C^2$$ is constant by the  construction of the solution, whereas
 $${\mathcal E}(V_{r_1}, E_{r_1})=(a  x + b )^2+ (c  x + d )^2$$ has a minimum
 $ x_* (t) = - \frac {a  b  + c  d } {2 (a^2 + c^2 )} \in (x_- (t), x_+ (t)) $.
 Since $ {\mathcal E}(V_{r_1}, E_{r_1})=  C^2  $ at $x=x_\pm (t)$,
 then
 $$ {\mathcal E}(V_{r_1}, E_{r_1})< C^2={\mathcal E}(V_{r_2}, E_{r_2}), \quad x(t)\in (x_- (t), x_+ (t)), \quad t\in (0, T_*).$$

 Both principles, on the basis of which admissible solutions can be distinguished, lead to the same conclusion:
 for the complete system \eqref {K1},
 the solution $ (V_{r_1}, E_{r_1}) $ must be chosen as a rarefaction wave, while when the condition \eqref {KE} is applied, only the possibility $ (V_ {r_2}, E_ {r_2}) $ remains.


 \section{Construction of a singular shock wave}\label{Shock}

We need to build a shock wave for the second part of the period $2 \pi$, $ t
 \in (T_*, 2 \pi) $. However, in order not to complicate the notation, we, without loss of generality, shift the time point $T_*$ to zero. Thus, we are in a situation where the initial data correspond to  a shock wave and $t=T_*$ is the point of the first intersection of the characteristics.

 Suppose that for $ t\in (0,T_*) $ we have constructed a solution to the problem as
 \begin{equation}\label{K15}
 (V_{s}, E_{s})=\left\{\begin{array}{cc}V_-(t), &x<{\Phi}(t),\\
  V_+(t), &x>{\Phi}(t), \end{array}\right.
\end{equation}
that is, we found the position of the shock wave $ x = {\Phi} (t) $. Then the density can be found as $n(t,x)=1-[E]|_{x={\Phi}(t)} \delta(x-{\Phi}(t))$.

Thus, we must take into account the presence of a strongly singular component of the solution.
However, before proceeding to the construction of a solution in this case, we will give a general definition of a strongly singular solution and obtain its main properties.

\subsection{Definition of a generalized strongly singular solution}

 Starting from the divergent form \eqref {K4}, we define a generalized strongly singular solution to the problem \eqref {K4}, \eqref {K2} according to  \cite{Shelkv}.

 Let
 \begin{eqnarray}\label{K19}
   V(t,x)&=&V_-(t,x)+[V(t,x)]|_{x={\Phi}(t)} \Theta(x-{\Phi}(t)), \\  E(t,x)&=&E_-(t,x)+[E(t,x)]|_{x={\Phi}(t)} \Theta(x-{\Phi}(t)),\label{K20}
   \\ n(t,x)&=&\hat n(t,x)+e(t)\delta(x-{\Phi}(t)),\label{K21}
    \end{eqnarray}
     where $[f]=f_+-f_-$, $f_\pm $ are differentiable functions having one-sided limits, $t\ge 0$, $x\in \mathbb R$, $\hat n(t,x)=1-\{E_x(t,x)\}$, $\{E_x\}$ is the derivative of the function $ E $ at the points at which it exists in the usual sense, $e(t):=e(t,{\Phi}(t))$,  $e(t)=-[E(t,x)]|_{x={\Phi}(t)}$.
\bigskip

\begin{definition}
The triple of distributions $ (V, E, n) $, given as \eqref {K19} - \eqref {K21} and the curve $ \gamma $, given as $ x = \Phi (t), $ $ \Phi (0) = 0 $, $ \Phi (t) \in C ^ 1 $, is called a generalized singular solution of the problem \eqref{K4},
\begin{eqnarray*}\label{K30}
(V,E,n)|_{t=0}=\\(V^0_-(x)+[V(x)]^0 \Theta (x),\, E^0_-(x) +[E(x)]^0 \Theta(x),\, n^0(x)=\hat n^0(x)+e^0\delta(x)),\nonumber
\end{eqnarray*}
if for all test functions $\phi(t,x)\in \mathcal{D}({\mathbb R}\times [0,\infty))$
\begin{eqnarray}\nonumber
\int\limits_0^\infty\int\limits_{\mathbb R} \hat n  (\phi_t+V \phi_x) dx dt +\int\limits_{\gamma} e(t) \frac{\delta \phi (t,x)}{\delta t} \frac{dl}{\sqrt{1+(\dot\Phi(t))^2}}+\\
\int\limits_{\mathbb R} \hat n^0(x) \phi(0,x) dx  + e(0) \phi(0,0)  =0,\label{K190}\\
\int\limits_0^\infty\int\limits_{\mathbb R} \left(( \frac{\hat n V^2}{2}+E^2)\phi_t+ \frac{\hat n V^3}{2} \phi_x\right) dx dt +\int\limits_{\gamma} \frac{e(t)(\dot \Phi(t))^2}{2} \frac{\delta \phi (t,x)}{\delta t}  \frac{dl}{\sqrt{1+(\dot\Phi(t))^2}}+\nonumber\\
\int\limits_{\mathbb R} \left(\frac{\hat n^0(x) (V^0(x))^2}{2} + (E^0(x))^2 \right)\phi(0,x) dx  + \frac{e(0)(\dot \Phi(0))^2}{2} \phi(0,0) =0,\label{K200}
\end{eqnarray}
 where $\int\limits_{\gamma} \cdot dl$ is the curvilinear integral along the curve $ \gamma $, the delta-derivative $\frac{\delta \phi (t,x)}{\delta t}\big|_{\gamma} $ is defined as the tangential derivative on the curve $ \gamma $, namely
$$
\frac{\delta \phi (t,x)}{\delta t}\big|_{\gamma} =\left(\frac{\partial \phi (t,x)}{\partial t}+ \dot\Phi(t) \frac{\partial \phi (t,x)}{\partial x}\big|_{\gamma} \right)\big|_{\gamma}=\frac{d \phi (t,\Phi(t))}{d t}= \sqrt{1+(\dot\Phi(t))^2} \frac{\partial \phi (t,x)}{d {\bf l}},
$$
where ${\bf l}=(-\nu_2, \nu_1)=\frac{(1, \dot\Phi(t))}{\sqrt{1+(\dot\Phi(t))^2}}$ is a unit vector tangent to $\gamma$.

\end{definition}
The action of the delta function $ \delta (\gamma) $ concentrated on the curve $ \gamma $ on the test function is defined according to \cite {Kanwal},
as $$
(\delta(\gamma),\phi(t,x))=\int\limits_{\gamma} \phi(t,x) \frac{dl}{\sqrt{1+(\dot\Phi(t))^2}},
$$
where $\phi(t, x)\in \mathcal{D}({\mathbb R}\times [0,\infty))$.

\bigskip

 \subsection{Rankine-Hugoniot conditions for delta-shock waves (the Rankine-Hugoniot deficit)}

 \begin{theorem} Let the domain $ \Omega \in {\mathbb R}^2 $ be divided by a smooth curve $ \gamma_t = \{(t, x): x = \Phi (t) \} $ into the left and right sides $ \Omega_\mp $. Let the triple of distributions $ (V, E, n) $, given as \eqref {K19} - \eqref {K21} and the curve $ \gamma_t $ be a strongly singular generalized solution for the system \eqref {K4}. Then this solution satisfies the following analogue of the Rankine-Hugoniot conditions
 \begin{eqnarray}\label{RH1}
   \frac{d}{dt} e(t)&=&\left(-[ \hat n  V]+[\hat n] \dot\Phi(t)\right)\big|_{x=\Phi(t)},\\
   \label{RH2}
  \frac{d}{dt}\frac{ e(t) (\dot \Phi(t))^2}{2}&=&\left(-\left[\frac{\hat n  V^3}{2}\right]+\left[\frac{\hat n V^2+E^2}{2}\right] \dot\Phi(t)\right)\big|_{x=\Phi(t)}.
 \end{eqnarray}
 \end{theorem}

The proof of the first statement, \eqref {RH1}, is contained in \cite {Shelkv}, the proof of \eqref {RH2} repeats the proof of the analogue of the Rankine-Hugoniot conditions for the energy equation in the "pressureless" gas dynamics model \cite {NRS}. Let us briefly recall this.

We denote  ${\bf n}=( \nu_1, \nu_2)=\frac{(\dot\Phi(t),-1)}{\sqrt{1+(\dot\Phi(t))^2}}$ the unit normal to the curve $ \gamma_t $ directed from $\Omega_-$ to $\Omega_+$.

Choose a test function $ \phi (t, x) $ with support $ K \subset \Omega $. Then
\begin{eqnarray*}
&&\int\limits_0^\infty\int\limits_{\mathbb R} \left(( \frac{\hat n V^2}{2}+E^2)\phi_t+ \frac{\hat n V^3}{2} \phi_x\right) dx dt =\\\nonumber
&&\int\limits\limits_{\Omega_-\cap K}\left(( \frac{\hat n V^2}{2}+E^2)\phi_t+ \frac{\hat n V^3}{2} \phi_x\right) dx dt +
\int\limits\limits_{\Omega_+\cap K} \left(( \frac{\hat n V^2}{2}+E^2)\phi_t+ \frac{\hat n V^3}{2} \phi_x\right) dx dt.
\end{eqnarray*}
Integration by parts by the second equation \eqref {K4} gives
\begin{eqnarray}\nonumber
&&\int\limits_{\Omega_\pm\cap K} \left(( \frac{\hat n V^2}{2}+E^2)\phi_t+ \frac{\hat n V^3}{2} \phi_x\right) dx dt =
-\int\limits_{\Omega_\pm\cap K} \left(( \frac{\hat n V^2}{2}+E^2)_t +( \frac{\hat n V^3}{2})_x\right) \phi \,dx dt\mp\\
\nonumber&&
\int\limits_{\gamma_t} \left(\nu_2 ( \frac{\hat n_\pm (V_\pm)^2}{2}+(E_\pm)^2) +\nu_1 \frac{\hat n_\pm(V_\pm)^3}{2}\right) \phi(t,x) dl-
\int\limits_{\Omega_\pm\cap K\cap \mathbb{R}} \left( \frac{\hat n^0(x)(V^0(x))^2}{2}+(E^0)^2\right) \phi(0,x) dx.
\end{eqnarray}
Thus,
\begin{eqnarray}\label{del3}
&&\int\limits_0^\infty\int\limits_{\mathbb R} \left(( \frac{\hat n V^2}{2}+E^2)\phi_t+ \frac{\hat n V^3}{2} \phi_x\right) dx dt +\\
\nonumber
&&\int\limits_{\Omega_\pm\cap K\cap \mathbb{R}} \left( \frac{\hat n^0(x)(V^0(x))^2}{2}+(E^0)^2\right) \phi(0,x) dx=\\\nonumber&&
 - \int\limits_{\gamma_t} \left( \left[ \frac{\hat n V^2}{2}+E^2\right]\nu_2 +\left[ \frac{\hat n V^3}{2}\right]\nu_1  \right) \phi(t,x) dl.
\end{eqnarray}

Further,
\begin{eqnarray}\label{del4}
&&\int\limits_{\gamma} \frac{e(t)(\dot \Phi(t))^2}{2} \frac{\delta \phi (t,x)}{\delta t}   \frac{dl}{\sqrt{1+(\dot\Phi(t))^2}}=\\
\nonumber
&&- \int\limits_{\gamma} \frac{\delta }{\delta t} \left( \frac{e(t)(\dot \Phi(t))^2}{2} \right) \phi(t,x) \frac{dl}{\sqrt{1+(\dot\Phi(t))^2}}
 - \frac{e(0)(\dot \Phi(0))^2}{2} \phi(0,0).
\end{eqnarray}
Adding \eqref {del3} and \eqref {del4}, taking into account \eqref {K200}, we get
\begin{eqnarray}\nonumber
\int\limits_{\gamma_t} \left( \left[ \frac{\hat n V^2}{2}+E^2\right]\nu_2 +\left[ \frac{\hat n V^3}{2}\right]\nu_1 -
  \frac{\delta }{\delta t} \left( \frac{e(t)(\dot \Phi(t))^2}{2} \right)\,\frac{1}{\sqrt{1+(\dot\Phi(t))^2}}\right)
 \phi(t,x) dl=0
\end{eqnarray}
for any $ \phi \in \mathcal D (\Omega) $. This implies \eqref {RH2}.

We see that the generalized Rankine-Hugoniot conditions are a system of second-order ordinary differential equations, therefore, to solve the Cauchy problem \eqref {K1}, \eqref {K1} (with the divergent form \eqref {K4}), \eqref {K19} , \eqref {K20} should set the initial velocity of the shock position $\dot \Phi(0)$.

Since the system \eqref {K1} has coinciding eigenvalues
$ \lambda_1 (V) = \lambda_2 (V) = V $, the admissibility condition for a singular shock wave coincides with the geometric entropy condition:
\begin{eqnarray}
\label{accept}
 \min\{V_-, V_+\} \le  \dot \Phi(t) \le \max\{V_-, V_+\},
\end{eqnarray}
meaning that characteristics from both sides come to the shock.

As we will see below, this condition allows us to obtain a condition on the derivative at an intermediate point and construct a solution to the Riemann problem in a unique way. In addition, in this problem, a final point arises, where the trajectory of the delta-shaped singularity must come, which also determines the problem.

\subsection{Mass and energy transfer ratios for a singular shock}
Suppose that $ V, E $ is a compactly supported classical solution of the \eqref {K1} system. Then, according to \eqref {K4}, the total mass  $$\int\limits_{\mathbb R}\,n(t,x) dx$$ and the total energy    $$\frac12\int\limits_{\mathbb R}\,(n(t,x)V^2(t,x) +  E^2(t,x))\, dx$$ are conserved.  Note that the total energy consists of kinetic and potential parts. Let us show that in order to obtain analogs of these conservation laws for a strongly singular solution, it is necessary to introduce the mass and energy concentrated on the shock. Suppose that the line of discontinuity $ x = \Phi (t) $ is a smooth curve.

We denote
\begin{eqnarray*}
 \label{M}
 &&\mathcal M(t)=\int\limits_{-\infty}^{\Phi(t)}\,n(t,x) dx +\int\limits_{\Phi(t)}^{+\infty} \,n(t,x) dx,
  \\
  \label{Ek}
&&  {\mathcal E}_{k}(t)=\frac12\,\left(\int\limits_{-\infty}^{\Phi(t)}\,{n(t,x)V^2(t,x)}\,dx +\int\limits_{\Phi(t)}^{+\infty} \,{n(t,x)V^2(t,x)}\, dx\right),\\
  \label{Eр}
&&  {\mathcal E}_{p}(t)=\frac12\,\int\limits_{-\infty}^{+\infty} E^2(t,x)\,dx,
\end{eqnarray*}
the mass, kinetic and potential energies concentrated outside the shock.
We interpret the amplitude $ e (t) $ and the term $ \frac12 {e (t) (\dot \Phi (t))^2} $ as the mass $ m (t) $ and kinetic energy $ w (t) $ concentrated on the shock.

\begin{theorem}
Let the solution \eqref {K19} - \eqref {K21} be a strongly singular solution to the system \eqref {K4}. Then the following
relations of balance take place:
\begin{eqnarray}
\label{balance}
\dot m(t)=-\dot {\mathcal M}(t), \quad  \dot w(t)=-\dot {\mathcal E}(t),\\
\label{conserv}
M(t)+m(t)=\mathcal M(0)+m(0), \quad \mathcal E(t)+w(t)=\mathcal E(0)+w(0).
\end{eqnarray}
\end{theorem}

Proof. Both equalities \eqref {balance}  can be proved in the same way. Let us prove, for example, the first of them. Because
\begin{eqnarray}
\nonumber
&&  \dot {\mathcal M}(t)=
   -[\hat n] \big|_{x=\Phi(t)}\dot\Phi(t)+
  \left(\int\limits_{-\infty}^{\Phi(t)}+\int\limits_{\Phi(t)}^{+\infty}\right)\,n_t(t,x) dx=\\
  \nonumber
 && - [\hat n] \big|_{x=\Phi(t)}\dot\Phi(t)-
  \left(\int\limits_{-\infty}^{\Phi(t)}+\int\limits_{\Phi(t)}^{+\infty}\right)\,(n(t,x)V(t,x))_x dx =\\
  \nonumber
&&  - [\hat n] \big|_{x=\Phi(t)}\dot\Phi(t)+
  [ \hat n V]\big|_{x=\Phi(t)},
\end{eqnarray}
together with \eqref {RH1} this equality shows that
\begin{eqnarray*}
  \dot m(t)+  \dot {\mathcal M}(t)=0,
\end{eqnarray*}
whence the first equalities in \eqref {balance} and \eqref {conserv} follow.

\bigskip

\subsection{Constructing a strongly singular shock wave for piecewise constant initial data \eqref{K3}}

We proceed to constructing a strongly singular solution in our particular case.

Since in this situation $ \hat n = 1 $ and, accordingly, $ [\hat n] = 0 $, the equations according to which it is possible to determine the amplitude and location of a strongly singular shock wave are greatly simplified and take the form
\begin{eqnarray}\label{RH11}
   \dot e(t)&=&-[V]\big|_{x=\Phi(t)},\\
   \label{RH21}
  \frac{d}{dt}{ e(t) (\dot \Phi(t))^2}&=&\left(-\left[{V^3}\right]+\left[{V^2+E^2}\right] \dot\Phi(t)\right)\big|_{x=\Phi(t)}.
 \end{eqnarray}
 Since the values of $ V_\pm (t) $ and $ E_\pm (t) $ are known (see \eqref {KvE}), the  values of jumps can be calculated directly:
  \begin{eqnarray*}\label{jumps}
  && [V^3]\big|_{x=\Phi(t)}=
    (([E]^0)^3-3 [E V^2]^0) \cos^2 t \sin t  \\
    &&\nonumber + (( [V]^0)^3-3 [E^2 V]^0)   \cos^3 t
   + 3 [E V^2]^0 \cos t - ([E]^0)^3 \sin t,
   \end{eqnarray*}
   and
  \begin{eqnarray*}
   \left[{V^2+E^2}\right]\big|_{x=\Phi(t)}=\left[V^2+E^2\right]^0\big|_{x=x_0}=K=\rm const.
\end{eqnarray*}
Therefore, from \eqref {RH1} we find
\begin{eqnarray}\label{ee}
    e(t)=-[V]^0 \sin t-[E]^0 \cos t.
 \end{eqnarray}
Note that for all $ t $ for which a shock wave exists, for $ e (0)> 0 $ the amplitude $ e (t) $ remains positive.

On the interval $ (0, T_*) $ there is  a point $ t_*$ at which $ V_-(t_*) = V_ + (t_*): = U $. Then from the admissibility condition \eqref{accept} we have
$ V_-(t_*) = \dot \Phi (t_*) = V_+(t_*)=U $.
It can be readily found that
 $$
 T_*= 2 \arctan\frac{[V]^0}{[E]^0},\quad t_*=\frac{T_*}{2},\quad U=
 \frac{E_-^0 V_+^0-E_+^0 V_-^0}{\sqrt{([V]^0)^2+([E]^0)^2}}.
 $$

We denote $(\dot \Phi (t)=q$ and $ (\dot \Phi (t))^2=Q(t)$.
 From \eqref{RH21} we get the Cauchy problems
 \begin{eqnarray}\label{q}
  \dot q&=&\frac{ - \left[{V^3}\right] +K q -\dot e q^2}{2 e q}, \quad q(t_*)=U,
 \end{eqnarray}
 and
\begin{eqnarray}\label{Q}
  \dot Q&=&-\frac{\left[{V^3}\right]}{e}+{\rm sign}\,U \,  \frac{K}{e}\sqrt{Q} -\frac{\dot e}{e}Q, \quad Q(t_*)=U^2.
 \end{eqnarray}

 The solutions to problems \eqref{q}, \eqref{Q} cannot be found explicitly, however, they always exist for $q\ne 0$.
 If  at a certain point $t=t_0\in (0,T_*)$, we have $q\to 0$, then, as follows from \eqref{q},\eqref{Q}, $\dot q\to \infty$, $\dot Q\to c={\rm const}\ne 0$. Indeed, if $c=0$, then $\left[{V^3}\right]|_{t=t_0}=0$. Nevertheless, it is easy to check that  $\left[{V^3}\right]=0$ if and only if $t=t_*$. This implies that $q=O(\sqrt{|t-t_0|})$, as $t\to t_0\ne t_*$.

Thus, if there exists a point $t=t_0\in (0,T_*)$ such that $Q(t_0)=0$, we first find the unique solution to \eqref{Q} at the first segment $(0,t_0)$ or $(t_0,T_*)$ (the segment must contain $t_*$), and then find the unique solution to the Cauchy problem
\begin{eqnarray*}\label{Q1}
\dot Q&=&-\frac{\left[{V^3}\right]}{e}-{\rm sign}\,U \frac{K}{e}\sqrt{Q} -\frac{\dot e}{e}Q, \quad Q(t_*)=0,
 \end{eqnarray*}
 on the second segment. Then we find $q$ on both segments.

 Let us note if $\Phi (t)$ has an extremum on $(0, T_*)$, then $\Phi(t)\in C^1(0,T_*)$ and can be found uniquely, however, $\ddot{\Phi} (t_0)$ does not exist.

 \subsection{Examples}

1. We start with the case \eqref{KE}, when the system \eqref{K5} can be reduced to one equation and one of the possible conservative form
is
\begin{equation*}\label{Kcons}
V_t+(\frac{V^2}{2})_x=-\sigma\sqrt{C^2-V^2},\quad \sigma={\rm sign} (-V_x)=\pm 1,
\end{equation*}
which does not require an introduction of a singular shock. The position of a usual shock is defined by the Rankine-Hugoniot condition and gives
\begin{equation}\label{RH11}
\dot \Phi(t)=\frac{V_-(t)+V_+(t)}{2}.
\end{equation}

Let us choose the initial data as
\begin{equation*}\label{RD1}
V_-^0=1,\,V_+^0=0,\, E_-^0=0,\,E_+^0=-1.
\end{equation*}
Here $T_*=\frac{\pi}{2}$, $K=U=0$ and $e(0)=1$.

\begin{center}

\begin{figure}[htb]

\hspace{-1cm}
\begin{minipage}{0.3\columnwidth}
\includegraphics[scale=0.25]{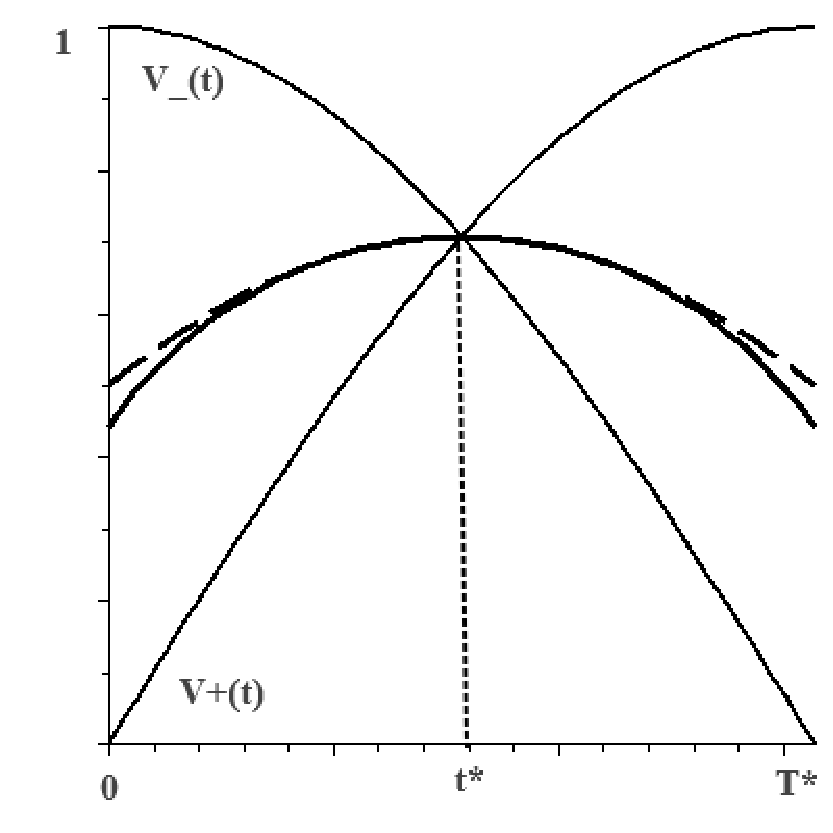}
\end{minipage}
\begin{minipage}{0.3\columnwidth}
\includegraphics[scale=0.25]{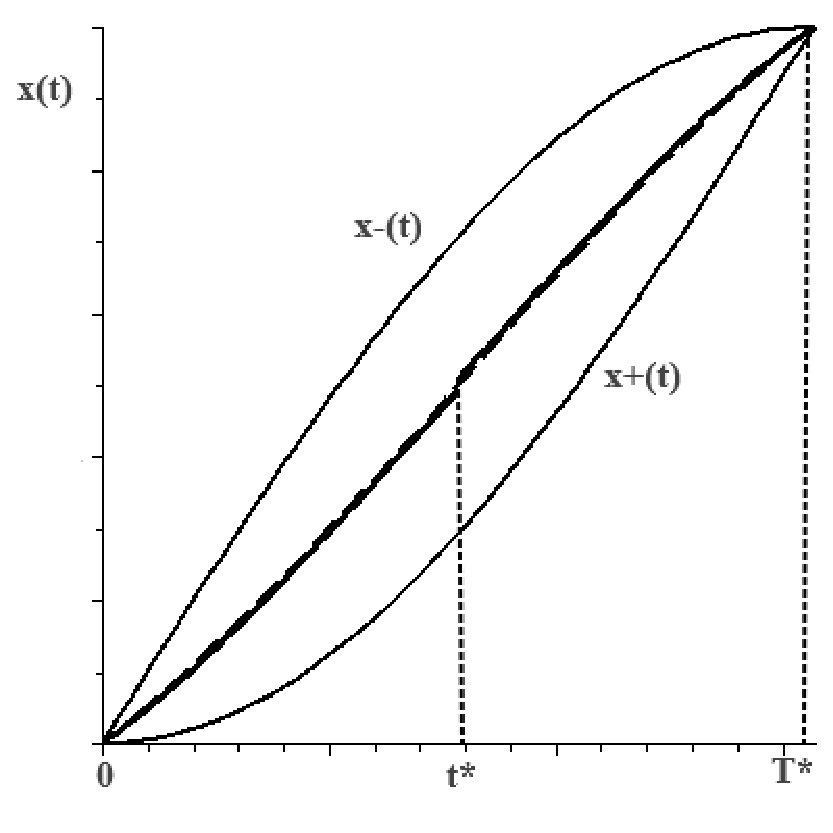}
\end{minipage}
\begin{minipage}{0.3\columnwidth}
\includegraphics[scale=0.25]{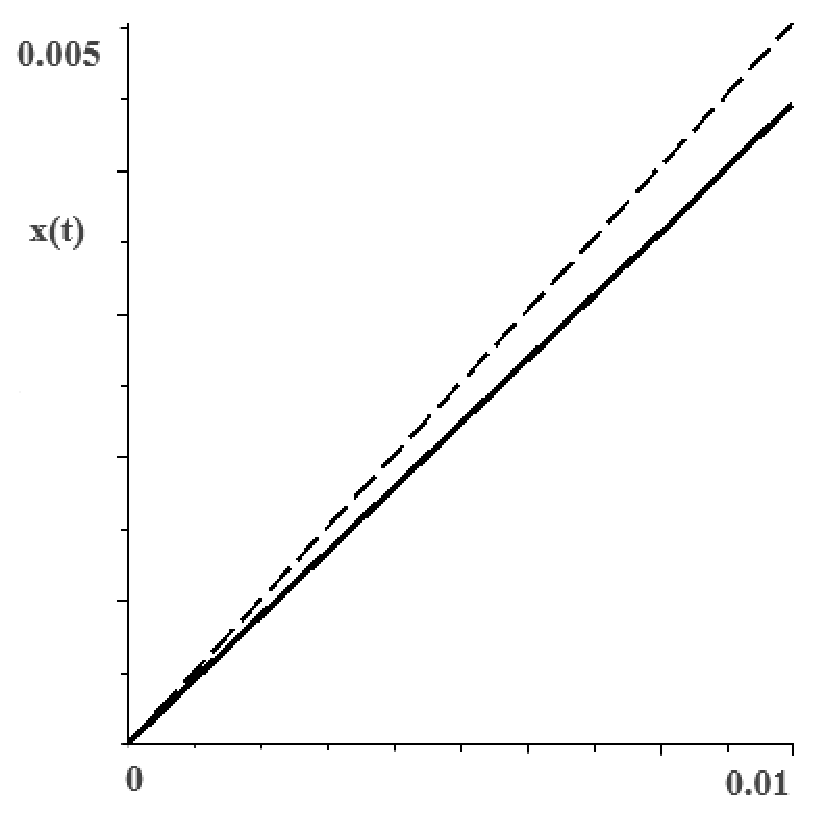}
\end{minipage}
\caption{Velocity (left) and position (center and right) of the singular shock (solid) vs. the velocity and position of the Rankine-Huhoniot shock (dash) for Example 1.}\label{Pic3}
\end{figure}
\end{center}

Fig.\ref{Pic3}, left, presents the behavior of the velocity $\dot \Phi(t)$ of the singular shock satisfying the geometric entropy condition \eqref{accept} (solid), in comparison with the velocity of shock based on the Rankine-Hugoniot condition \eqref{RH11} (dash). One can see that the difference is very small. Fig.\ref{Pic3}, center, presents the position of the singular shock
between characteristics (solid), in comparison with the Rankine-Hugoniot shock (dash), the difference is almost negligible. Fig.\ref{Pic3}, right, shows the zoom of this difference near the origin.

\begin{center}

\begin{figure}[htb]

\hspace{-1cm}
\hspace{1cm}
\begin{minipage}{0.4\columnwidth}
\includegraphics[scale=0.25]{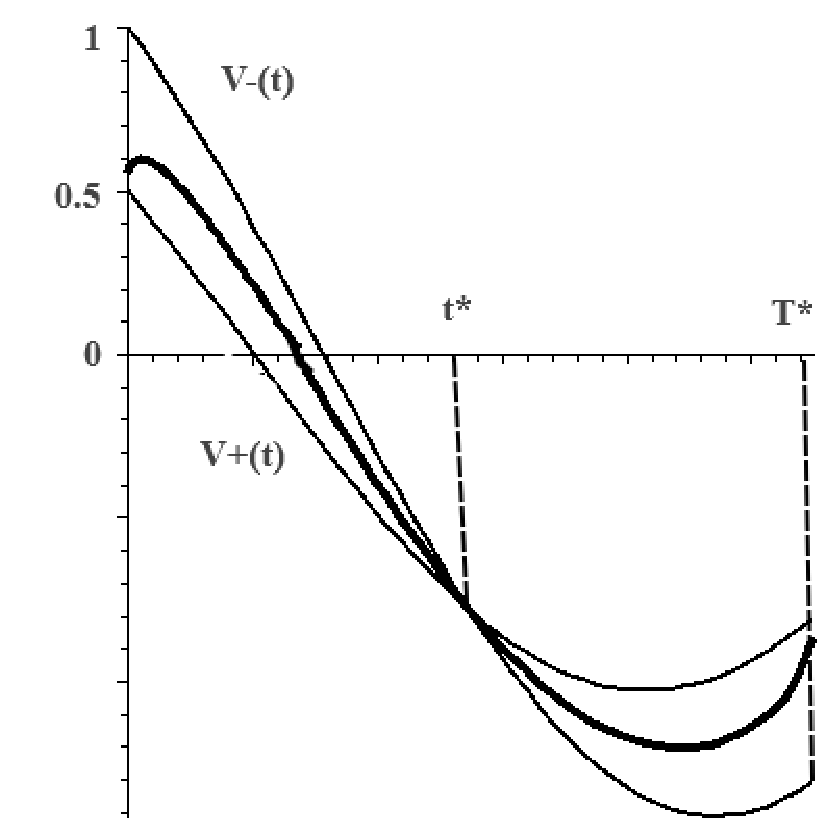}
\end{minipage}
\hspace{1.5cm}
\begin{minipage}{0.4\columnwidth}
\includegraphics[scale=0.25]{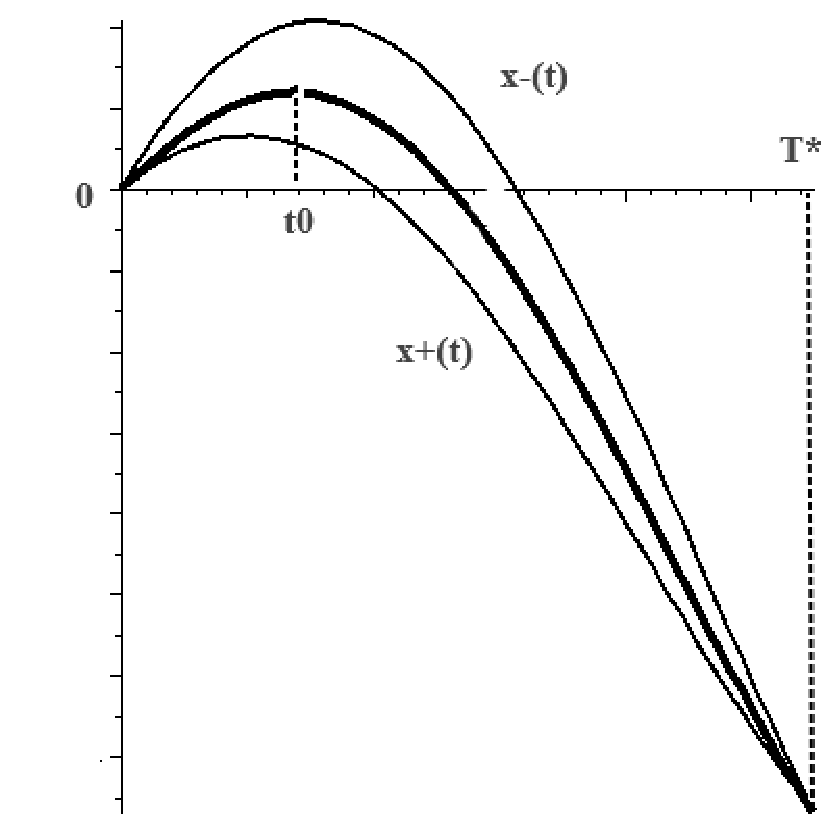}
\end{minipage}
\caption{ Velocity (left) and position (right) of the singular shock for Example 2. }\label{Pic4}
\end{figure}
\end{center}

2.  The next example is for the following data:
\begin{equation*}\label{RD2}
V_-^0=1,\,V_+^0=0.5,\, E_-^0=1,\,E_+^0=0.9.
\end{equation*}
Here $T_*=2.746801534$, $U=-.7844645404$, $K=-.94$ and $e(0)=0.1$. This example is interesting, since $\dot \Phi(t)$ changes the sign at a point $t_0=0.69174927\ne t_*$.
Fig.\ref{Pic4}, left, presents the behavior of the velocity $\dot \Phi(t)$ of the singular shock satisfying the geometric entropy condition \eqref{accept}, Fig.\ref{Pic2}, right, presents the position of the singular shock
between characteristics.

The solution in the examples are found numerically by means of the Runge-Kutta-Fehlberg method of fourth-fifth order.

 \section{Discussion}\label{Dis}

1. We show that the reduced equations of a cold plasma provide the simplest example of an inhomogeneous system in which the solution of the Riemann problem consists of a rarefaction wave and a shock wave periodically replacing each other. The system is not so interesting from a physical point of view, since it is believed that the cold plasma equations are valid only for a smooth solution \cite{Ch_book}. However, they are extremely interesting mathematically. Indeed, firstly, due to the non-strict hyperbolicity of the system, one can construct an example of non-uniqueness of the rarefaction wave. Second, the natural conservative form of the system can be used to construct a singular shock wave.

2. The solution of the Cauchy problem \eqref{K1}, \eqref{K3} can be rewritten in terms of the solution of the Euler-Poisson equation \eqref{EP} with discontinuous data $(n_0, V_0)$.

3. A similar procedure for solving the Riemann problem can also be applied to other non-strict hyperbolic systems written initially in a non-divergent form. The method consists in introducing an ``artificial density'', which makes it possible to write the system in a conservative form and define a strong singular solution.

4. The appearance of multiple rarefaction waves was noticed earlier in other models, for example, \cite{MM}.

5. The presence of pressure apparently prevents the existence of a strong singular solution \cite{Haspot}, in other words, the situation is similar to the influence of pressure in the gas dynamics model without pressure.

6. It should be noted that there are different plasma models and a huge amount of literature devoted to shock waves there. One of the popular models is the VlasovЦMaxwell system, which describes a collisionless ionized plasma \cite{Morawetz}, \cite{Schaeffer}, \cite{Polovin}, \cite{Balog}. Another assumption about plasma naturally changes the properties of shock waves.

7. The non-strictly hyperbolic system considered here has a simple wave solution  (an invariant manifold),  what makes them related to systems of the Temple class \cite{Temple}. However, the most interesting feature of the solution of the Riemann problem for
  the system of cold plasma equations are rooted in its inhomogeneity, while the equations of the Temple class are homogeneous and have constant states to the left and right of the shock wave or rarefaction wave.

 \section*{Acknowledgements} Supported by  the Moscow Center for Fundamental and Applied Mathematics under the agreement є075-15-2019-1621. The author thanks her former student Darya Kapridova for numerical calculations confirming the construction of a rarefaction wave in the case of a simple wave solution. The author expresses her sincere gratitude to the anonymous referee for careful reading.

\end{document}